\documentclass{colen2e}

\begin{document}

\def\[[#1]]{[\![#1]\!]}
\def\((#1)){(\!(#1)\!)}
\def\olm{\leq_{\mathrm{lex}}}
\def\olM{\geq_{\mathrm{lex}}}
\def\tv#1{#1_1,#1_2,#1_3}
\def\vg#1#2{(#1_1,\ldots,#1_{#2})}
\def\fio{for\-ma\-lly in\-de\-pen\-dent o\-ver}
\def\sg{sys\-tem of ge\-ne\-ra\-tors}
\let\em=\sl

\title{Valuations in fields of power series}

\author{Herrera Govantes, Francisco Javier\thanks{Partially supported
    by FQM--218 (JdA).}\\
Departamento de \'Algebra \\ Universidad de Sevilla\\ E-mail:
jherrera@algebra.us.es \and Miguel \'Angel Olalla
Acosta$^{*\ddag}$\thanks{Partially supported by HF2000-0044.}\\
Departamento de \'Algebra \\ Universidad de Sevilla\\ E-mail:
olalla@algebra.us.es \and Jos\'e Luis Vicente
C\'ordoba$^{*}$\thanks{Partially supported by BFM
    2001--3207 (MCyT).}\\
Departamento de \'Algebra \\ Universidad de Sevilla\\ E-mail:
jlvc@us.es }

\date{September, 2001}
\maketitle

\thispagestyle{empty}

\begin{abstract}
This paper deals with valuations of fields of formal
meromorphic functions and their residue fields. We
explicitly describe the residue fields of the monomial
valuations. We also classify all the discrete rank one
valuations of fields of power series in two and three
variables, according to their residue fields. We prove
that all our cases are possible and give explicit constructions.
\end{abstract}

\keywords{Valuation Theory, Power Series}

\mathclass{12J20}{13A18}

\section{Introduction}

In  this  paper,  we  give  a  ``classification'' of certain
valuations  of  $k\((X_1,\ldots,X_n))$,  where  $k$  is  an
algebraically  closed  field,  namely  discrete   valuations
finite over $k\[[X_1,\ldots,X_n]]$  and having as  center in it
the maximal ideal.

In section 2, we construct explicitly the residue field of the
discrete monomial valuations of any rank. Section 3 is devoted to
list in detail all the discrete rank one valuations of $k\((X))$,
$k\((X_1,X_2))$, $k\((X_1,X_2,X_3))$ (again with the condition on
the center). The case of $k\((X))$ is trivial.

In the case  of $k\((X_1,X_2))$ we get that the residue field
$\Delta_v$ of any such valuation $v$ is a pure transcendental
extension $k(u)$ of $k$ generated by one element, and $v$ itself
is in one of the following two cases:
\begin{enumerate}

\item either $v$ is monomial, or

\item $v$ is the composition of a finite number of
blowing-ups and coordinate changes with a monomial valuation.

\end{enumerate}

In the case of $k\((X_1,X_2,X_3))$, the situation is much more
complicated. To compute the residue field and to give an explicit
expression of $v$, we allow ourselves to perform sequences (maybe
infinite) of blowing-ups and coordinate changes. The possible
cases are the following:
\begin{enumerate}

\item $\Delta_v=k(u)$. In this case, $k\[[X_1,X_2,X_3]]$ can
be embedded into a power series ring in two variables contained in
the ring of the completion of $v$, and the extension to it of $v$
is monomial.

\item $\Delta_v$ has transcendence degree $2$ over $k$, and
the algebraic part may be non-trivial. This case includes the
monomial valuations and the compositions of a finite sequence of
blowing-ups with such a valuation.

\item $\Delta_v$ has transcendence degree $1$ and the
algebraic part may be finite or countably infinite. If $u$ is
transcendental over $k$,  every countably generated algebraic
extension of $k(u)$ can be realized as the residue field
$\Delta_v$ of the valuation $v$.

\end{enumerate}

Section 4 is devoted to discrete rank two valuations of
$k((X_1,X_2))$. The principles of the process of threatament are
similar to the rank one case. The result is that any discrete rank
two valuation of $k((X_1,X_2))$ is zero--dimensional and:

\begin{enumerate}
\item either it's monomial,
\item or it can be reduced to a monomial one by a (possibly infinite)
  sequence of blowing-ups and coordinate changes.
\end{enumerate}

\section{Monomial valuations}

Throughout all this paper, $k$ will be an algebraically closed
field, $R=k\[[\mathbf{X}]]=k\[[X_1,\ldots,X_n]]$ the formal power
series ring in the variables $X_i$ and  $K$ its quotient field.

\begin{remark}{\label{11}}

{\bf 1)} Every $f\in R$ will be written as
$f=\sum_{A\in\mathbf{Z}_{0}^{n}}f_A \mathbf{X}^A$, where, if
$A=(a_1,\ldots,a_n)$, then $\mathbf{X}^A$ means $X_1^{a_{1}}\cdots
X_n^{a_{n}}$. We will also write $$
\mathcal{E}(f)=\{A\in\mathbf{Z}_{0}^n\mid f_A\neq 0\} \, . $$

{\bf 2)} We will freely speak about valuations of $R$, meaning the
restriction to $R$ of a valuation of its quotient field. To
abbreviate, the word ``valuation'' will replace the phrase
``discrete $k$-valuation of $R$ centered at the maximal ideal
$M=({\bf X})\cdot R$''.

{\bf 3)} We will  use in every $\mathbf{Z}^m$ the lexicographic
order, which will be denoted by $\olm$. It is a monomial order in
the sense that it is compatible with the additive group structure.

{\bf 4)} Let $0<m\leq n$ be an integer  and $$
L=\{B_1,\ldots,B_n\}\subset \mathbf{Z}_0^m\setminus\{\mathbf{0}\}
$$ a \sg\ of $\mathbf{Z}^m$. We associate to each monomial
$\mathbf{X}^A$ in $R$ an element of $\mathbf{Z}_0^m$, which we
call its $L$-{\em degree}\index{degree}, defined in the following
way: $$ \mathrm{degree}_{L}(\mathbf{X}^A)=\sum_{i=1}^n a_iB_i\, ,
\quad A=(a_1,\ldots,a_n)\, . $$
\end{remark}

\begin{lemma}{\label{11a}} Let
$B=(b_1,\ldots,b_m)\in\mathbf{Z}_0^m$; then the set of monomials
with $L$-degree equal $B$ is finite. Consequently, the homogeneous
components of the $L$-{\em graduation}\index{graduation} of $R$
are finite dimensional  $k$-vector spaces.
\end{lemma}

\begin{proof}
It is enough to observe that, since all the $B_i$ are different
from zero, any linear combination $B=\sum_{i=1}^n a_iB_i$ with
$a_i\in\mathbf{Z}_0$ must have the coeficients bounded by the
maximum of the $b_j$.
\end{proof}

\begin{definition}{\label{12}} Let $0<m\leq n$ be an integer
and $$ L=\{B_1,\ldots,B_n\}\subset
\mathbf{Z}_0^m\setminus\{\mathbf{0}\} $$ a  \sg\ of
$\mathbf{Z}^m$. Let $v:\,R\to\mathbf{Z}^m\cup\{\infty\}$ be the
function which assigns $\infty$ to zero and $$
v(f)=\min_{\mathrm{lex}}\{\mathrm{degree}_{L}(\mathbf{X}^A) \mid
A\in \mathcal{E}(f)\} $$ to $0\neq f\in R$. The extension of $v$
to $K$  is a valuation of $K/k$ whose value group is
$\mathbf{Z}^m$, which will be called the monomial
valuation\index{monomial valuation} associated to $L$.
\end{definition}

\begin{remark}{\label{14}}
Let $A_1,\ldots ,A_p\in\mathbf{Z}^n$; the following conditions are
equivalent:
\begin{enumerate}
\item[a)] The monomials $\{\mathbf{X}^{A_1},\ldots ,\mathbf{X}^{A_p}\}$ are
algebraically independent over $k$ (they may have negative
exponents).
\item[b)] $\{A_1,\ldots ,A_p\}$ are {\bf Z}--linearly independent
in $\mathbf{Z}^n$.
\end{enumerate}
\end{remark}

\begin{proposicion}{\label{15}}
Let $v$ be a valuation of $K/k$ of rank $m$. Let $B_i=v(X_i)$ for
$i=1,\ldots ,n$ and consider the set $L=\{ B_1,\ldots , B_n\}$; we
will further assume that $L$ is a \sg\ of $\mathbf{Z}^m$. Let
$v_L$ be the corresponding monomial valuation of $K/k$. Then, for
every $f\in R$, one has $v(f)\olM v_L(f)$.
\end{proposicion}

\begin{proof} Let us assume that $f\neq 0$ and let
$v_L(f)=B$. Let $\mathbf{X}^A$ be a monomial with non-zero
coefficient in $f$ of $L$-degree $B$ and consider the fraction
$g/\mathbf{X}^A$. It is clear that $g/\mathbf{X}^A\in R[Y_1,\ldots
,Y_s]$, where each $Y_i$ is a quotient of monomials in
$\mathbf{X}$ of value zero. It is also clear that $R[Y_1,\ldots
,Y_s]\subset R_v$, hence $v(g/\mathbf{X}^A)\olM 0$, so $v(g)\olM
v(\mathbf{X}^A)=B$. This completes our proof.
\end{proof}

\begin{theorem}{\label{16}} Let $v$ be a valuation of
$K/k$ of rank $m$,  $B_i=v(X_i)$ for $i=1,\ldots ,n$ and let us
assume that $L=\{ B_1,\ldots, B_n\}$ is a \sg\ of $\mathbf{Z}^m$.
Then the following conditions are equivalent:
\begin{enumerate}
\item $v$ is the monomial valuation associated to $L$.
\item The residue field $\Delta_v$ of $v$ is the
extension of $k$ (pure transcendental if $n>m$) generated by the
$n-m$ monomials (fractional in general) whose exponents constitute
a basis of the solutions of the system $\sum_{i=1}^n
y_iB_i=\mathbf{0}$, as an abelian group.
\end{enumerate}
\end{theorem}

\begin{proof}
Here we are using the known fact (although not completely trivial)
that the solutions in $\mathbf{Z}^m$ of a system of homogeneous
linear equations is a free abelian group (c.f. \cite{jlv}). Let us
assume that $v$ is the monomial valuation associated to $L$ and
let $K_L$ be the subfield of $K$ consisting of the quotients of
$L$-forms of the same $L$-degree. The natural map $$
\begin{array}{ccc}
K_L & \rightarrow  &\Delta_v \\ \displaystyle{\frac{f_B}{g_B}} &
\rightarrow &\displaystyle{\frac{f_B}{g_B}+{\cal M}_v}
\end{array}
$$ is injective, so we may assume $K_L\subset\Delta_v$. Let
$f/g\in K$ be such that $v(f/g)=0$, that is, $f$  and
 $g$ are two power series such that $v(f)=v(g)=B$.

Let us write $f=f_B+f_1$ and $g=g_B+g_1$, where $f_B$ and $g_B$
are their $L$-initial forms respectively. Then, $$
\frac{f}{g}-\frac{f_B}{g_B}=\frac{g_Bf_1-f_Bg_1}{g_Bg}\, , $$ so
$v(f/g-f_B/g_B)>0$ and $$ \frac{f}{g}+{\cal
M}_v=\frac{f_B}{g_B}+{\cal M}_v\, , $$ therefore $K_L=\Delta_v$.

Let $M$ be the submodule of $\mathbf{Z}^m$ consisting of the
solutions of the system $\sum_{i=1}^n y_iB_i=\mathbf{0}$ and
$\{A_1,\ldots,A_{n-m}\}$ be a base of $M$. Then one easily checks
that $$ K_L=k(\mathbf{X}^{A_1},\ldots ,\mathbf{X}^{A_{n-m}})\, .
$$

Now let us assume that the residue field of $v$ is
$\Delta_v=k(\mathbf{X}^{A_1},\ldots ,\mathbf{X}^{A_{n-m}})$, where
$A_1,\ldots ,A_{n-m}$ form a basis of the submodule $M$ of the
solutions of the system $\sum_{i=1}^n y_iB_i=\mathbf{0}$. Let
$v_L$ be the monomial valuation associated to $L$. By proposition
\ref{15}, in order to show that $v=v_L$ it is enough to prove that
$v(f)=v_L(f)$ for every $L$-form $f$.

Assume the contrary and let $f$ be an $L$-form such that
$v(f)>_{\mathrm{lex}}v_L(f)$. Let us write $$ f=\sum_{i=1}^s
\alpha_i\mathbf{X}^{C_i} $$ where $\alpha_i\in k$ and the sum is
extended to all the monomials of $L$-degree $v_L(f)$. For every
$j=1,\ldots ,n$ one has $v(f/\mathbf{X}^{C_j})>0$, hence $$
\frac{f}{\mathbf{X}^{C_j}}+ {\cal M}_v=0+{\cal M}_v\, . $$ Since
$$ \frac{f}{\mathbf{X}^{C_j}}=\alpha_j+\sum_{i\neq
j}\alpha_i\mathbf{X}^{C_i-C_j} $$ one has $\sum_{i\neq
j}\alpha_i\mathbf{X}^{C_i-C_j}+{\cal M}_v= -a_j+{\cal M}_v$. This
implies that $\sum_{i\neq j}\alpha_i\mathbf{X}^{C_i-C_j}+{\cal
M}_v$ belongs to $k$, which is not possible since all the
$\mathbf{X}^{C_j-C_i}$ are monomials in the $A_l$ and these are
algebraically independent over $k$.
\end{proof}

\section{Discrete rank one valuations in low dimension}

We now consider other discrete valuations of fields of power
series. It is easy to make a complete list of all the discrete
rank one valuations of $k\((X_1))/k$  and $k\((X_1,X_2))/k$, and
we do it in a few considerations. However, the case of three
variables is more difficult. We will explicitly describe all the
discrete rank one valuations of $k\((\tv{X}))/k$ because it is the
most difficult one.

\begin{remark}{\label{21}}

{\bf 1)} As before, let $R=k\[[\mathbf{X}]]=k\[[X_1,\ldots,X_n]]$,
$n>1$, and  $K$ be its quotient field. We fix a discrete rank one
valuation $v:K=k\((X_1,\ldots,X_n))\to\mathbf{Z}\cup\infty$,
centered at the maximal ideal of $R$ and we assume that the value
group is $\mathbf{Z}$. We will denote, as usual, by $R_v$, ${\cal
M}_v$ and $\Delta_v$ the ring, the ideal and the residue field of
$v$.

{\bf 2)} We will consider the completion $\widehat{v}$ of $v$,
together with its ring $\widehat{R}_{\widehat{v}}$ and the
quotient field $\widehat{K}_{\widehat{v}}$ of
$\widehat{R}_{\widehat{v}}$. We will fix a datum, which will play
a key role in our study, namely a section
$\sigma:\Delta_v\to\widehat{R}_{\widehat{v}}$ of the natural
homomorphism $\widehat{R}_{\widehat{v}}\to\Delta_v$, which exists
by the Cohen structure theorem. {\em We will always identify
$\Delta_v$ with its image in $\widehat{R}_{\widehat{v}}$ by
$\sigma$, so we will assume from now on that
$\Delta_v\subset\widehat{R}_{\widehat{v}}$}.

{\bf 3)} Remark that, if $t\in\widehat{R}_{\widehat{v}}$ is an
element of value $1$, then $t$ is formally independent over
$\Delta_v$ and $\widehat{R}_{\widehat{v}}=\Delta_v\[[t]]$.

\end{remark}

\begin{remark}{\label{22}} Let us assume that $n=1$; then
the usual order function $\nu_{X_1}$  is a discrete rank one
valuation whose ring is $k\[[X_1]]$ and its residue field is $k$.
Every other discrete valuation $v$  of $k\((X_1))/k$ such that
$R_v\supset k\[[X_1]]$  and ${\cal M}_v\cap k\[[X_1]]=(X_1)$ must
coincide with $\nu_{X_1}$. Thus the only  non trivial valuation of
rank $m=1$ of $k\((X_1))/k$ is the usual order function.
\end{remark}

\begin{remark}{\label{24}} Let us assume that $v$ is a
discrete rank one valuation of $K/k$, where
$K=k\((X_1,\ldots,X_n))$, and refer to the above notations.

{\bf 1)} Let $z_1,z_2\in \widehat{{\cal M}}_{\widehat{v}}$ and $L$
a subfield of $\widehat{R}_{\widehat{v}}$ containing $k$, i.e.,
such that $k\subset L\subset\Delta_v$. Let us assume that
$\widehat{v}(z_2)>\widehat{v}(z_1)>0$; then, the natural map
$k\[[z_1,z_2/z_1]]\to \widehat{R}_{\widehat{v}}$ is injective, as
we will see in a moment (c.f. {\bf 3)}). We will assume that
$k\[[z_1,z_2/z_1]]\to \widehat{R}_{\widehat{v}}$. Therefore,
$z_2/z_1\in\widehat{{\cal M}}_{\widehat{v}}$, every power series
$f(z_1,z_2/z_1)\in L\[[z_1,z_2/z_1]]$ makes sense in
$\widehat{R}_{\widehat{v}}$,
$L\[[z_1,z_2/z_1]]\subset\widehat{R}_{\widehat{v}}$ and
$L\[[z_1,z_2]]\subset L\[[z_1,z_2/z_1]]$. This $L\[[z_1,z_2/z_1]]$
is called the blowing-up of $L\[[z_1,z_2]]$. If
$\widehat{v}(z_2)=q\widehat{v}(z_1)+r$ is the euclidean division
and $r>0$, then $q$ repetitions of the blowing-up dividing by
$z_1$ give us the ring $L\[[z_1,z_2/z_1^q]]\supset L\[[z_1,z_2]]$
in which $\widehat{v}(z_2/z_1^q)=r$. If the remainder is zero, we
usually take $q-1$ blowing-ups instead of $q$ just to equate the
values. The pair $(z_1,z_2/z_1^q)$ is monomial, birrational with
respect to $(z_1,z_2)$, i.e., each one is a monomial in the other,
possibly with negative exponents. Equivalently, the vectors of the
exponents $\{(1,0),(-q,1)\}$ form a basis of $\mathbf{Z}^2$.

{\bf 2)} Let $d=\gcd(v(X_1),\ldots,v(X_n))$; then a finite
sequence of bowing-ups (describing the Euclid's algorithms to
compute the greatest common divisors in the sense of {\bf 1)} will
produce a vector $\mathbf{z}=\vg{z}{n}$ of elements in
$\widehat{{\cal M}}_{\widehat{v}}$ such that $\widehat{v}(z_i)=d$,
$\forall i=1,\ldots,n$. The important point here is that
$k\[[\mathbf{z}]]\supset k\[[\mathbf{X}]]$ and the vectors
$\mathbf{z}$ and $\mathbf{X}$  are monomial birrational.


This is obviously true for every starting vector
$\mathbf{x}=\vg{x}{n}$ with components in $\widehat{{\cal
M}}_{\widehat{v}}$, either variables (i.e., formally independent)
or not, over any starting field. This process will be called {\em
reduction of a vector to the minimum value}.

{\bf 3)} In the process of reduction of a vector to the minimum
value, a crucial point is that, if the starting vector has two
components formally independent over the ground field $L$, then
the two components of the final vector are also formally
independent over $L$. To see this, it is enough to prove that, if
$z_1,z_2$ are formally independent over $L$, then $z_1,z_2/z_1$
are also formally independent over $L$. Reasoning by
contradiction, let us assume that $z_1,z_2/z_1$ are not formally
independent over $L$ and let $f(Z_1,Z'_2)\in L\[[Z_1,Z'_2]]$ be
such that $f(z_1,z_2/z_1)=0$; then $f$ cannot be a unit. By the
Weierstra\ss\ preparation theorem,  there is a unit $u(Z_1,Z'_2)$
and a non-unit Weierstra\ss\ polynomial $g(Z_1,Z'_2)$ in $Z'_2$
such that $f(Z_1,Z'_2)=Z_1^r u(Z_1,Z'_2)g(Z_1,Z'_2)$. Thus $$
0=f(z_1,z_2/z_1)=z_1^r u(z_1,z_2/z_1)g(z_1,z_2/z_1)\, , $$ so
$g(z_1,z_2/z_1)=0$. If $d=\mathrm{degree}_{Z'_2}\,(g)$, then
$0=z_1^d g(z_1,z_2/z_1)=g'(z_1,z_2)$, which is not possible
because $z_1,z_2$ are formally independent and $g'(z_1,z_2)$ is a
monic polynomial in $z_2$.
\end{remark}

\begin{remark}{\label{25}} The above remarks  \ref{24} allow
us to describe completely the discrete rank one valuations of
$k\((X_1,X_2))/k$ with our initial conditions.

{\bf 1)} We start with $(X_1,X_2)$ and reduce it to the minimum
value, obtaining a vector $\mathbf{z}=(z_1,z_2)$ whose components
are formally independent over $k$ by remark \ref{24}.{\bf 3)}, and
$k\[[z_1,z_2]]\supset k\[[X_1,X_2]]$. If the residue of $z_2/z_1$
belongs to $k$, we denote it by $\alpha$ and take the element
$z'_2=z_2-\alpha z_1$ whose value is strictly greater that
$\widehat{v}(z_1)$. The vector  $(z_1,z'_2)$ has components again
formally independent over $k$ and $k\[[z_1,z'_2]]=k\[[z_1,z_2]]$.

{\bf 2)} If $\widehat{v}(z'_2)$ is a multiple of
$\widehat{v}(z_1)$ we again reduce to the minimum value, which is
$\widehat{v}(z_1)$, taking blowing-ups dividing by $z_1$. If,
again, the residue is rational, we repeat and so on. This process
cannot be infinite because it would amount to an expansion of
$z_2$ as a power series in $z_1$ with coefficients in $k$, which
is not possible by formal independence (c.f. remark \ref{24}.{\bf
3)}). Therefore, the process stops, either because we arrive at an
element whose value is not a multiple of $\widehat{v}(z_1)$ or we
arrive at a residue which is transcendental over $k$. In the first
case, we reduce again to minimum value getting an element of value
strictly smaller than $\widehat{v}(z_1)$ and we start over from
the beginning. In the second case, we stop. This game of falling
in the first case can be repeated only  finite number of times,
because droppings of positive value can be only finitely many so,
after a finite number of steps, we get a vector,  renamed
$\mathbf{z}=(z_1,z_2)$, in which  the components have the same
value and the residue $z_2/z_1$ is transcendental over $k$.
Moreover, and this is the most important fact, these components
are formally independent over $k$ and $k\[[z_1,z_2]]\supset
k\[[X_1,X_2]]$. Let $u_2\in\widehat{R}_{\widehat{v}}$  be the
residue of $z_2/z_1$ and write  $z_2=u_2z_1+z'_2$ with
$\widehat{v}(z'_2)>\widehat{v}(z_1)$. Let us call $d$ the value
$d=\widehat{v}(z_1)$.

{\bf 3)} The restriction  $w$ of $\widehat{v}$ to $k\[[z_1,z_2]]$
is a monomial valuation. In fact, let $0\neq f(z_1,z_2)\in
k\[[z_1,z_2]]$ be a non-unit of order $r>0$ and let $$
f(z_1,z_2)=\sum_{i\geq r} f_i(z_1,z_2) $$ its decomposition into
sum of forms; then, if we consider the inclusion
$k\[[z_1,z_2]]\subset k(u_2)\[[z_1]]$, we can write $$
f(z_1,z_2)=f(z_1,u_2z_1+z'_2)=f_r(z_1,u_2z_1)+T=z_1^r f_r(1,u_2)+T
$$ where  $f_r(1,u_2)\neq 0$ because $u_2$ is transcendental over
$k$ and $T$ has value greater than $rw(z_1)=rd$ by proposition
\ref{15}. This proves our claim.

{\bf 4)} However, {\bf 3)} proves more. From it, we see that every
non-unit power series of $k\[[z_1,z_2]]$ has a value which is a
multiple of $d$. In particular, every non-unit power series in
$k\[[X_1,X_2]]$ has the same property. Since the value group of
$v$ is $\mathbf{Z}$, we must have $d=1$. We also see that the
residue of every element in $K\((X_1,X_2))$ is a rational function
of $u_2$, hence $\Delta_v=k(u_2)$. Finally $v$ coincides with the
composition of the (injective) ring homomorphism
$\varphi:k\[[z_1,z_2]]\to k(u_2)\[[t]]$ given by $\varphi(z_1)=t$,
$\varphi(z_2)=u_2t$ with the $t$-order function on $k(u_2)\[[t]]$.
\end{remark}

These remarks prove the following theorem, due to Briales and
Herrera \cite{Br-He,Bri2}.

\begin{theorem}{\label{25k}}
A discrete rank one valuation $v$ of $k\[[X_1,X_2]]$ either is
monomial or it can be reduced to a monomial one by a finite
process consisting of blowing-ups and coordinate changes.
\end{theorem}

The rest of this section is devoted to listing all possible
discrete rank one valuations of $k\((\tv{X}))/k$ whose ring
contains $k\[[\tv{X}]]$ and its center here is the maximal ideal.
We fix one of them, $v$ and proceed.

\begin{remark}{\label{26}}

{\bf 1)} We start with the vector $\mathbf{X}=(\tv{X})$ and reduce
it to its minimum value $d$; we get a vector $\mathbf{z}=(\tv{z})$
with $\widehat{v}(z_i)=d$, $i=1,2,3$. Let $\alpha_i$ be the
residue of $z_i/z_1$, $i=2,3$; if both are elements of $k$, we
take $z'_i=z_i-\alpha_i z_1$ and the vector
$\mathbf{z}'=(z_1,z'_2,z'_3)$. If $\widehat{v}(z'_2)$ and
$\widehat{v}(z'_3)$ are multiples of $d$, we reduce $\mathbf{z}'$
to its minimum value by blowing-up dividing by $z_1$; let us
rename $\mathbf{z}=(\tv{z})$ the output vector. If, again, the two
residues belong to $k$ and the new values are multiple of $d$ we
continue. The important point here is that, always,
$k\[[\mathbf{X}]]\subset k\[[\mathbf{z}]]$.

{\bf 2)} Is it possible to enter in an infinite process of this
kind? In other words, is it possible to arrive to
$k\[[\mathbf{X}]]\subset k\[[\mathbf{z}]]=k\[[z_1]]$? The answer
is no. In fact, if affirmative, it will imply that all the $X_i$
belong to $k\[[z_1]]$ so, by Weierstrass preparation arguments,
for each $i=2,3$ the $X_i$ will satisfy an equation of integral
dependence over $k\[[X_1]]$, which is not the case.

{\bf 3)} Let us assume that, after a finite number of steps, we
get a vector $\mathbf{z}$ such that, either $\widehat{v}(z_2)$ or
$\widehat{v}(z_3)$ is not a multiple of $\widehat{v}(z_1)$. Then,
reduction to minimum value will give us a vector, again renamed
$\mathbf{z}$, such that the common value of its three components
is $d'<d$. If, again, we get residues in $k$ as in {\bf 1)} and
enter into a cycle as in there, we see by {\bf 2)} that the cycle
cannot be infinite.

{\bf 4)} If, again, the minimum value drops and we enter into a
cycle as in {\bf 1)}, and so on, we see that this process cannot
be infinite, either. The reason is that a decreasing sequence of
positive integers must stabilize.

{\bf 5)} Therefore, after a finite number of steps, we must arrive
at a vector $\mathbf{z}$ of equally valued components such that,
after reordering if needed, the residue of $z_2/z_1$ is
transcendental over $k$. In this process, we could have ``lost''
the component $z_3$. For instance, let us consider two variables
$x,y$, the ring $k\[[x,y]]$ and, inside it, our starting ring
$k\[[\tv{z}]]$ with $z_1=y$, $z_2=xy$, $z_3=yf(x)$,where $f(x)$ is
transcendental over $k(x)$. These $z_1,z_2,z_3$ are formally
independent over $k$, as it is easy to see. Let us consider the
monomial valuation $v$ of $k\((\tv{z}))$ such that $v(z_1)=1$,
$v(z_2)=2$, $v(z_3)=1$; then $x=z_2/z_1$ has value $1$,
$k\[[z_1,z_2/z_1]]=k\[[x,y]]$ and $z_3\in k\[[x,y]]$. Therefore,
in all events, we could continue with the ring $k\[[z_1,z_2/z_1]]$
simply forgetting the $z_3$. In any case, deciding whether $z_3$
is lost or not could be the outcome of an infinite process. In
fact, to see that $z_3\in k\[[z_1,z_2]]$ or not involves a
(possibly infinite) process of blowing-ups.

{\bf 6)} Let us assume that, in addition to the assumptions of
remark \ref{26}.{\bf 5)}, $z_3\notin k\[[z_1,z_2]]$. Repeating a
process similar to the one in {\bf 5)}, this time only with
$z_1,z_3$, we could arrive either to a new fall of minimal value
or to a residue transcendental over $k$. In case of fall of
minimal value, we start everything from the beginning, and so on.
It is evident that this process must be finite.

{\bf 7)} The end of the history has two possibilities:
\begin{enumerate}
\item[a)] either a vector $\mathbf{z}=(z_1,z_2)$ of equally
valued components such that the residue of $z_2/z_1$ is
transcendental over $k$ and $$ k\[[\tv{X}]]\subset k\[[z_1,z_2]]\,
, $$
\item[b)] or not the above; then we have  a vector
$\mathbf{z}=(\tv{z})$ of equally positive valued components such
that the residues of $z_2/z_1$ and $z_3/z_1$ are transcendental
over $k$ (not necessarily algebraically independent over $k$) and
$$ k\[[\tv{X}]]\subset k\[[\tv{z}]]\, . $$
\end{enumerate}
\end{remark}

\begin{remark}{\label{27}}
In this note we consider the first case of the output in remark
\ref{26}.{\bf 7)}. We have that $k\[[\tv{X}]]\subset
k\[[z_1,z_2]]$ and the residue $u$ of $z_2/z_1$ is transcendental
over $k$. Let $d=\widehat{v}(z_1)=\widehat{v}(z_2)$ and write
$z'_2=z_2-uz_1\in\widehat{R}_{\widehat{v}}$ with
$\widehat{v}(z'_2)>d$. Let $\mathbf{Z}=(Z_1,Z_2)$ be a vector of
indeterminates and $0\neq f\in k\[[Z_1,Z_2]]$ be a power series of
order $r>0$ and denote by $f_i(Z_1,Z_2)$ the forms of $f$; then $$
f(z_1,z_2)=\sum_{i\geq r} f_i(z_1,uz_1+z'_2)=f_r(z_1,uz_1)+T
=z_1^rf_r(1,u)+T\, , $$ where $f_r(1,u)\neq 0$ by transcendence of
$u$  and $T\in\widehat{R}_{\widehat{v}}$, $\widehat{v}(T)>r$. This
shows that $\widehat{v}(f(z_1,z_2))=rd$ and, since
$k\[[\tv{X}]]\subset k\[[z_1,z_2]]$ and the value group of $v$ is
$\mathbf{Z}$, then it must be $d=1$. On the other hand, this also
proves that the residue field of the restriction $w$ of
$\widehat{v}$ to $k\[[z_1,z_2]]$ is $k(u)$, which is $\Delta_v$,
is equal to $k(u)$. Of course, this restriction $w$ is a monomial
valuation.

We can always write a concrete example of such a valuation by
taking two formally independent variables $z_1,z_2$ a monomial
valuation of $k\[[z_1,z_2]]/k$ and three formally independent
power series in $k\[[z_1,z_2]]/k$.
\end{remark}

\begin{remark}{\label{28}}
In these notes we consider the second case of the output in remark
\ref{26}.{\bf 7)}. We have that $k\[[\tv{X}]]\subset k\[[\tv{z}]]$
and the residues $u_{2,1}$ of $z_2/z_1$ and $u_{3,1}$ of $z_3/z_1$
are transcendental over $k$. In this case we need a little more
preparation. Let $d=\widehat{v}(z_1)$.

{\bf 1)} Now we initiate a process of coordinate changes and
blowing-ups, similar to the one in remarks \ref{26}, in search of
an element of value strictly smaller than $d$, if it exists. We
start with $(z_1,z_2)$, take the extension $L_1=k(u_{2,1})$ and
the element $z'_2=z_2-u_{2,1}z_1$. If $\widehat{v}(z'_2)$ is a
multiple of $\widehat{v}(z_1)$, then we equate values by taking
blowing-ups dividing by $z_1$. Then we take again a suitable
extension $L_2=L_1(u_{2,i_{2}})$ and perform a coordinate change
$z''_2=z_2-u_{2,1}z_1-u_{2,i_{2}}z_1^{i_{2}}$. If
$\widehat{v}(z''_2)$ is a not a multiple of $\widehat{v}(z_1)$,
then we take blowing-ups dividing by $z_1$ until we get an element
of strictly smaller value. In this case, we restart everything
(reduction to minimum value, coordinate changes, and so on,
starting from a new three components vector), again over $k$ as
the ground field.

For $(z_1,z_2)$, if we fall into an infinite process  of values
multiples of $\widehat{v}(z_1)$,  we have a power series expansion
$$ z_2=\sum_{j\geq 1} u_{2,j}z_1^j\, , \quad u_{2,j}\in
\Delta_v\setminus\{0\} $$ and we act likewise with $(z_1,z_3)$.

{\bf 2)} This time, the end of the history is two power series
expansions $$ z_i=\sum_{j\geq 1} u_{i,j}z_1^j\, , \quad i=2,3\,
\quad u_{i,j}\in \Delta_v\setminus \{0\} $$ where, of course, we
have renamed the vector $\mathbf{z}$ and $u_{2,1},u_{3,1}$ are
transcendental over $k$. We denote by $L$ the extension of $k$
generated by all the $u_{i,j}$. If we write $z_2=u_{2,1}z'_1$ then
the map sending $z_1$ onto $z'_1$ and $z_i$ onto $z_i$, $i=2,3$ is
an $L$-automorphism of $L\[[\mathbf{z}]]$, so we may assume that
$z_2=uz_1$, $u=u_{2,1}$. Therefore, we have
\begin{equation}{\label{ec1}}
\begin{array}{ccl}
z_2 & = & uz_1 \\ z_3 & = & \displaystyle{\sum_{j\geq 1}
u_{3,j}z_1^j}
\end{array}
\end{equation}
so $L\[[\tv{z}]]\subset L\[[z_1]]\subset\widehat{R}_{\widehat{v}}$
and $L\subset\Delta_v$.

{\bf 3)} In the situation of {\bf 2)}, we consider the restriction
$w$ of $\widehat{v}$ to $L\[[z_1]]$; then $w$ is necessarily the
$z_1$-order function, so $\widehat{v}(z_1)=1$. Moreover,
$\Delta_v=L$ and $\widehat{R}_{\widehat{v}}=L\[[z_1]]$.
\end{remark}

So, in these remarks we have proved the following theorem:

\begin{theorem}{\label{28k}}
Let $v$ be a discrete rank one valuation of $k\((\tv{X}))/k$,whose
ring contains $k\[[\tv{X}]]$, its center here is the maximal ideal
and the group of values is ${\bf Z}$. Then we have one of the
following situations:
\begin{enumerate}
\item[A)] There exists a vector ${\bf z}=(z_1,z_2)$ of elements in
  $\widehat{\cal M}_{\widehat{v}}$ such that $k\[[\tv{X}]]\subset
  k\[[z_1,z_2]]$, $v(z_1)=v(z_2)=1$, $\Delta_v=k(u)$, where $u$ is the
  residue of $z_2/z_1$, and the restriction $w$ of $\widehat{v}$ to
  $k\(({\bf z}))$ is a monomial valuation.
\item[B)] There exists a vector ${\bf z}=(z_1,z_2,z_3)$ of elements in
  $\widehat{\cal M}_{\widehat{v}}$ such that $k\[[\tv{X}]]\subset
  k\[[z_1,z_2,z_3]]$, $v(z_1)=v(z_2)=v(z_3)=1$, with $z_2=uz_1$ and
  $z_3=\sum_{j\geq 1}u_{3,j}z_1^j$, where $u$ and $u_{3,1}$ are
  transcendental residues over $k$, $\Delta_v=k(u,\{ u_{3,j}\}_{j\geq
  1})$ and $\widehat{R}_{\widehat{v}}=\Delta_v\[[z_1]]$.
\end{enumerate}
In both cases the vector ${\bf z}$ is explicitly obtained from
${\bf X}$ by a process consisting of blowing-ups and coordinate
changes.
\end{theorem}

\begin{remark}{\label{29}}
In the situation of equation (\ref{ec1}) the naturally arising
problem is to study the field extension $k(u)\subset\Delta_v$. In
this remark we deal with the case in which this extension is
transcendental.

{\bf 1)} By assumption, one of the coefficients $u_{3,j}$ must be
transcendental over $k(u)$; let us call $u'$ the first $u_{3,j}$
which is transcendental over $k(u)$, that is the one with the
smallest possible $j=j_0$. If there are $u_{3,j}\neq 0$  with
$j<j_0$ then we take the finite algebraic extension
$L=k(u)(\{u_{3,j}\}_{1\leq j<j_0})\supset k(u)$ and write
\begin{eqnarray*}
z'_3 & = & \displaystyle{\sum_{j=1}^{j_0-1}u_{3,j} z_1^j}\\ z''_3
& = & u'z_1^{j_0}+\displaystyle{\sum_{j>j_{0}}u_{3,j} z_1^j}
\end{eqnarray*}
We get the element $(1/z_1^{j_{0}-1})z''_3$ after a finite
sequence of change of coordinates and blowing-ups dividing by
$z_1$.

{\bf 2)} Let us assume that $z'_3$ does not exist, i.e.,
$z_3=u'z_1^j+\cdots$ ; then the components of the vector
$\mathbf{z}=(\tv{z})$ are formally independent over $k$ and
$\Delta_v=k(u,u')$. To see this, let $\mathbf{Z}=(\tv{Z})$ be a
vector of variables, $0\neq f\in k\[[\mathbf{Z}]]$ be a non-unit
of order $r$ and let $f_i(\mathbf{Z})$ be the form of degree $i$
of $f$; then $$ f(\tv{z})=\sum_{i\geq r} f_i(\tv{z})= z_1^r
f_r(1,u,u')+T $$ where $T\in\widehat{{\cal
M}}_{\widehat{v}}^{r+1}$ and $f_r(1,u,u')\neq 0$. This proves that
$f(\tv{z})\neq 0$ and its residue is a polynomial in $u,u'$, which
implies our claim.

{\bf 3)} Let us assume that $z'_3\neq 0$ and let $L'\supset L$
(c.f. remark \ref{29}.{\bf 1)}) be the minimal Galois extension of
$k(u)$ containing $L$. Let $$ P(z_1,Z_3) =\prod
(Z_3-(z'_3)^{(i)})\,  $$ be the product taken over the different
conjugates of $z'_3$. Then $P(z_1,Z_3)$ is a polynomial in the
indeterminate $Z_3$ with coefficients in $k(u)[z_1]$ such that
$P(z_1,z'_3)=0$ and, if $Q(z_1,Z_3)\in k(u)[z_1]$ is such that
$Q(z_1,z'_3)=0$ then $P(z_1,Z_3)|Q(z_1,Z_3)$.

Now, let $0\neq f(\tv{Z})\in k\[[\tv{Z}]]$ be an irreducible
non-unit such that
\begin{eqnarray*}
0 & = & f(\tv{z})=f(z_1,uz_1,z'_3+z''_3) \\ & = &
\displaystyle{\sum_{i\geq 0} f_i(z_1,uz_1,z'_3) (z''_3)^i}\, .
\end{eqnarray*}
This implies that $$ f_i(z_1,uz_1,z'_3)=0\, , \quad \forall i\geq
0\, , $$ so $P(z_1,Z_3)|f_i(z_1,uz_1,Z_3)$ for all $i$ by the
Weierstrass preparation theorem. By irreducibility,
$f(z_1,uz_1,Z_3)$ must be a unit factor of $P$. Hence, the initial
form cannot be vanished by replacing $Z_3$ by $z'_3+z''_3$.
Therefore, $(\tv{z})$ are formally independent over $k$.

{\bf 4)} It is obvious that there exist such valuations, for
formally independent initial arguments. The composition of the
following substitutions with the $t$-order functions give
valuations, the first one being monomial and the second one being
not:
\begin{eqnarray*}
\varphi(X_1) & = & t\; , \varphi(X_2) =  ut \; , \varphi(X_3) =
u't \\ \varphi'(X_1) & = & t\; , \varphi'(X_2) =  ut \; ,
\varphi'(X_3) =  (\sqrt{u})t+u't^2
\end{eqnarray*}

\end{remark}

\begin{remark}{\label{210}} We end with the case in which
$\Delta_v$ is an algebraic extension of $k(u)$.

{\bf 1)} In remark \ref{26}.{\bf 5)} we saw how $\Delta_v$ can be
a finite algebraic extension of $k(u)$ but not on the initial
variables. It is also possible to get a non-trivial extension. Let
us consider two variables $x,y$, the ring $k\[[x,y]]$ and, inside
it, our starting ring $k\[[\tv{z}]]$ with $z_1=y^2$, $z_2=xy^2$,
$z_3=y\exp(x)$, which are formally independent over $k$, as it is
easy to see. Let us consider the monomial valuation $v$ of
$k\((\tv{z}))$ such that $v(z_1)=2$, $v(z_2)=3$, $v(z_3)=1$; then
$x=z_2/z_1$ has value $1$, $k\[[z_1,z_2/z_1]]=k\[[x,y^2]]$ and
$z_3$ satifies $z_3^2=z_1\exp{2(z_2/z_1)}$.

{\bf 2)} It is also possible that the extension $\Delta_v\supset
k(u)$ is infinite. In this case, $(\tv{z})$ are formally
independent over $k$ as we can easily prove,reasoning by
contradiction, using an Galois argument similar to the one in
remark \ref{29}.{\bf 3)}. In other words, every countably
generated algebraic extension of $k(u)$ can be realized as the
residue field $\Delta_v$ of the valuation $v$.
\end{remark}

\section{Discrete rank two valuations in dimension two}

The principles of the techniques we have employed so far can be
applied to other cases. We make a careful study of the discrete
rank two valuations of $k\[[X_1,X_2]]$, just to illustrate the
ideas. Let $v$ be a discrete rank two valuation of
$K=k\((X_1,X_2))$, finite over $R=k\[[X_1,X_2]]$ whose center at
$R$ is the maximal ideal. We assume that the value group of $v$ is
$\mathbf{Z}^2$, which is no restriction at all. This means, that
there exist $z_1,z_2\in K$ such that $v(z_1)=(1,0)$,
$v(z_2)=(0,1)$.

Remark that the rank of any discrete  valuation of $R$ is at most
$2$ because $2$ vectors must generate a submodule of maximal rank
of $\mathbf{Z}^2$. Remark, further, that the limitation of the
rank by the dimension is by no means a consequence of the rather
special starting situation. In fact, it is a determined by the
Abhyankar-Zariski inequality (c.f. \cite{Ab}), valid for any local
ring of finite dimension.

We observe that the valuation $v$ must be zero-dimensional. In
fact, if the transcendence degree of $\Delta_v/k$ were positive,
there could be a composite of the corresponding place with a
non-trivial place of the residue field, which would be a valuation
of higher rank. This is not possible by the limitation of the rank
by the number of variables.

\begin{remark}{\label{30}}
Let $u_1,u_2$ be variables, $F=k\((u_2))\((u_1))$,
$T=k\((u_2))\[[u_1]]$  and let us define a special valuation in
$F$ by giving its action on $T$.

{\bf 1)} Any element in $T$ can be written as 
$$ w=\sum_{i\geq 0}
w_i(u_2)u_1^i\, , \quad w_i(u_2)\in k\((u_2))\, . 
$$ 
We can write $w$ in another form as 
$$ 
w= \sum_{(i,j)\in\mathbf{Z}_0\times\mathbf{Z}} a_{ij}\, u_1^i\,
u_2^j\, , \quad a_{i,j}\in k 
$$ 
and we will also write 
$$
\mathcal{E}(w)=\{(i,j)\in\mathbf{Z}_0\times\mathbf{Z}\mid
a_{ij}\neq 0\} 
$$ 
Let us denote by $\nu_{u_1}$, $\nu_{u_2}$ the
usual order of a power series in one variable (which can be
negative). Assume $w\neq 0$; then it is very easy to see that 
$$
\min_{\mathrm{lex}}(\mathcal{E})=
(\nu_{u_1}(w),\nu_{u_2}(w_{\nu_{u_1}(w)}(u_2)))\, . 
$$ 
It also easy to see that  the map $$ 0\neq w\mapsto
\min_{\mathrm{lex}}(\mathcal{E}) $$ defines a discrete rank two
valuation $\hat{v}$ of $T$. Let us denote also by $\hat{v}$ the
extension to $F$ of this valuation. The ring $R_{\hat{v}}$
consists of all the power series in $F$ with non-negative order in
$u_1$ and positive or infinite order of the coefficient power
series of $u_1^0=1$.  Remark that $R_{\hat{v}}\subset T$ and
$R_{\hat{v}}\neq T$; in fact all the terms in $T$ with order zero
in $u_1$ and leading coefficient of negative order are out of
$R_{\hat{v}}$.

{\bf 2)} We can embed $R$ into $R_{\hat{v}}$ in such a way that
$v$ extends uniquely to $\hat{v}$. This $R_{\hat{v}}$ plays the
role of the completion in the rank one case.

{\bf 3)} Assume that $v(X_2)\geq v(X_1)$; then the embedding of
$R$ into $R_{\hat{v}}$ has a natural extension to an embedding of
the blowing-up ring $k\[[X1,X_2/X_1]]$ into $R_{\hat{v}}$. The
proof is similar to the one in remark \ref{24}.{\bf 3)}.
\end{remark}

\begin{remark}{\label{31}}

{\bf 1)}  Let  us  assume  that  $v(X_1)$   and $v(X_2)$ are
$\mathbf{Z}$-linearly independent; then each $L$-form with
$L=\{v(X_1),v(X_2)\}$  is a  monomial.   In particular,
$\{v(X_1),v(X_2)\}$ is a basis of $\mathbf{Z}^2$.  Therefore,
every such valuation  of rank  2 of  $k\((X_1,X_2))/k$ is monomial
and zero-dimensional.

{\bf 2)} Let  us  assume  that  $v(X_1)$ and $v(X_2)$ are
$\mathbf{Z}$-linearly dependent; then we apply a process of
reduction to minimum value and a change of coordinates, as above.
We repeat this again and again. After a possibly infinite sequence
of blowing-ups and coordinate changes, we fall in a new vector
$(y_1,y_2)$ in ${\cal M}_{\hat{v}}$ such that the values generate
$\mathbf{Z}^2$. Then the valuation is monomial.

{\bf 3)} We give an example of an infinite process. Let $v$ be the
composition of the embedding $$
\begin{array}{ccl}
R & \longrightarrow & R_{\hat{v}} \\ X_1 & \longmapsto & u_2
\\ X_2 & \longmapsto & [\exp(u_2)-1]+[\exp(u_1/u_2)-1]
\end{array}
$$ with $\hat{v}$. The term $[\exp(u_2)-1]$ makes it necessary to
take an infinite sequence of blowing-ups before finding a vector
$(y_1,y_2)$ as before.
\end{remark}

\begin{remark}{\label{32}} As a final remark, we point out
that the process of reducing to the minimum value is, in all
cases, similar to the exhaustion of the first characteristic
exponent in the local resolution of the singularity of an
analytically irreducible plane curve.
\end{remark}

\end{document}